\documentclass[11pt,leqno,a4paper]{article}

\usepackage{amsmath}
\usepackage{amssymb}
\usepackage{amsthm}
\usepackage{enumerate.mod}
\usepackage{ifthen}
\usepackage{calc}

\newtheorem{thm}{Theorem}[section]
\newtheorem{lem}[thm]{Lemma}

\newtheorem{claim}[thm]{Claim}

\theoremstyle{definition}
\newtheorem{defn}[thm]{Definition}

\newtheorem*{notn}{Notation}

\theoremstyle{remark}

\newcommand{\zfc}{\textrm{\rm ZFC}}
\newcommand{\zfcf}{\mathrm{ZFC}^*}

\newcommand{\Fin}{\operatorname{Fin}}

\newcommand{\gdot}{\star}
\newcommand{\power}{\P}
\newcommand{\lbrak}{\bigl\|}
\newcommand{\rbrak}{\bigr\|}
\newcommand{\lmeas}{\mu\bigl(\lbrak}
\newcommand{\rmeas}{\rbrak\bigr)}

\newcommand{\inv}{^{-1}}

\renewcommand{\div}{\mathbin{/}}

\newcommand{\diff}{\bigtriangleup}

\newcommand{\alg}{\operatorname{alg}}

\newcommand{\forces}{\mskip 5mu plus 5mu\|\hspace{-2.5pt}{\textstyle \frac{\hspace{4.5pt}}{\hspace{4.5pt}}}\mskip 5mu plus 5mu}

\newcommand{\Iff}{\espc\mathrm{iff}\espc}

\renewcommand{\And}{\espc\text{and}\espc}

\newcommand{\m}{\mathfrak m}

\newcommand{\two}{\{0,1\}}

\newcommand{\pN}{\power(\N)}

\newcommand{\N}{\mathbb N}

\newcommand{\random}{\R}

\newcommand{\ctbl}{{\C}{\T}}

\DeclareFontFamily{U}{cmsy}{}
\DeclareFontShape{U}{cmsy}{m}{n}{<12> sfixed * [10] cmsy10 
<10> <9> <8> <7> <6> <5> sfixed * [10] cmsy10}{}
\DeclareSymbolFont{customtwo}{U}{cmsy}{m}{n} 
\DeclareMathSymbol{\sctn}{\mathord}{customtwo}{"78}

\DeclareFontFamily{U}{cmmi}{}
\DeclareFontShape{U}{cmmi}{m}{n}{<20> sfixed * [11] cmmib10 <12> sfixed * [10]
cmmi10 <10> <9> <8> sfixed * [6] cmmi6 <5> <6> <7> sfixed * [5] cmmi5}{}
\DeclareSymbolFont{custom}{U}{cmmi}{m}{n}
\DeclareMathSymbol{\rharpoon}{\mathord}{custom}{"2A}
\newlength{\widt}
\newlength{\widttwo}
\newlength{\hgt}
\renewcommand{\vector}[1]{\settowidth{\widttwo}{$\rharpoon$}\addtolength{\widttwo}{-1.0pt}\settowidth{\widt}{$#1$}\settoheight{\hgt}{$#1$}{#1}\hspace{-\widt}\raisebox{\hgt}{$\rharpoon$}\addtolength{\widt}{-\widttwo}\hspace{\widt}}

\newcommand{\oone}{{\omega_1}}

\newcommand{\vlambda}{\varLambda}
\newcommand{\vgamma}{\varGamma}

\renewcommand{\>}{\rangle}
\newcommand{\wh}{\widehat}

\newcommand{\espc}{\quad}

\newcommand{\ulc}{\ulcorner}
\newcommand{\urc}{\urcorner}

\newcommand{\A}{\mathcal A}
\newcommand{\B}{\mathcal B}
\newcommand{\C}{\mathcal C}
\newcommand{\D}{\mathcal D}
\newcommand{\E}{\mathcal E}
\newcommand{\F}{\mathcal F}
\newcommand{\ideal}{\mathcal I}

\newcommand{\G}{\mathcal G}
\renewcommand{\H}{\mathcal H}

\renewcommand{\P}{\Pcal}
\newcommand{\Pcal}{\mathcal P}
\newcommand{\Q}{\mathcal Q}
\newcommand{\R}{\mathcal R}
\newcommand{\Scal}{\mathcal S}
\newcommand{\T}{\mathcal T}
\newcommand{\U}{\mathcal U}
\newcommand{\V}{\mathcal V}

\newcommand{\X}{\mathcal X}

\newcommand{\Rfrak}{\mathfrak R}

\newcounter{saveenumi}

\newcommand{\roone}{\random^{\oone}}

\newcommand{\pob}{(\ideal^+\times\ideal^+)\gdot\random}

\renewcommand{\oone}{\omega_1}

\newcommand{\stdrep}{\operatorname{\Rfrak}}

\renewcommand{\vector}{\vec}

\newenvironment{enumeq}[1][none]{
\ifthenelse{\equal{#1}{none}}
{\begin{enumerate}\setcounter{enumi}{\value{equation}}}
{\begin{enumerate}[#1]\setcounter{enumi}{\value{equation}}}}
{\setcounter{equation}{\value{enumi}}\end{enumerate}}

\begin{document}

\title{Some partition properties for measurable colourings of~$\oone^2$}
\author{James  Hirschorn\footnote{The author acknowledges the generous support of
    the Japanese Society for the Promotion of Science (JSPS Fellowship for Foreign
    Researchers, ID\# P04301).}}
\date{January 18, 2005}
\maketitle

\begin{abstract}
We construct a measure on $\oone^2$ over the ground model in the forcing extension of a measure algebra,
and investigate when measure theoretic properties of some measurable colouring
of $\oone^2$ imply the existence of an uncountable subset of $\oone$ whose square is
homogeneous. This gives a new proof of the fact that, under a suitable axiomatic
assumption, there are no Souslin $(\oone,\oone)$ gaps in the Boolean algebra $L^0(\mu)\div\Fin$ when
$\mu$ is a separable measure.
\end{abstract}

\section{Introduction}
\label{sec:introduction}

The purpose of this note---the title notwithstanding---is to begin the task of extending the classical
notion of a measure to the notion of a measure over some model of (some
fragment of) set theory. We demonstrate an application of this generalization of
measure theory. 

Let $\random$ denote a measure algebra. 
In $W^\random$ (the forcing extension of $W$ by~$\random$), 
where $W$ is a certain transitive extension of the ground model universe $V$, we construct a
measure over $V^\random$ on subsets of $\oone^V\times\oone^V$. We then investigate
some properties of subsets $K$ of $\oone^V\times\oone^V$, with respect to this measure,
which imply the existence of an uncountable subset $X\subseteq\oone$ in~$V^\random$
homogeneous for $K$. This gives a new proof
of the fact that $\m>\aleph_1$ implies the nonexistence of Souslin $(\oone,\oone)$
gaps in $L^0(\mu)\div\Fin$ for $\mu$ separable; in other words, the nonexistence of
a destructible $(\oone,\oone)$ gap in $\pN\div\Fin$ after adding one random real.

\subsection{Credits/Disclaimer}
\label{sec:creditsdiclaimer}

The author wishes to declare that Stevo Todor\v cevi\'c of University of Toronto,
University of Paris VII (and anyone else besides the author) 
made no contribution whatsoever to this research. In particular, this research was
never discussed with him in person, by e-mail, or by any other means, and any new
ideas presented here are due entirely to the author. 

\section{Measure Theory}
\label{sec:meastheory}

The reader may wish to consult~\cite{F3} as an additional reference for measure theory.
Let us recall the classical notion of a measure.
We follow the convention that a \emph{measure} on a set $X$ is a partial function $\mu$ from
a nonempty Boolean subalgebra $\A$ of $\power(X)$ into $[0,\infty]$ satisfying
\begin{enumerate}[(a)]
\item $\mu(\emptyset)=0$,
\item for every sequence $(E_n:n\in\N)$ of pairwise disjoint members of $\A$, if
  $\bigcup_{n=0}^\infty E_n\in\A$ then $\mu\bigl(\bigcup_{n=0}^\infty
  E_n\bigr)=\sum_{n=0}^\infty\mu(E_n)$,
\end{enumerate}
while a \emph{measure space} is a triple $(X,\B,\mu)$ where $\B$ is a \emph{$\sigma$-algebra}
of subsets of $X$ and $\mu:\B\to[0,\infty]$ is a measure on $X$. 

Generalizing beyond algebras of subsets, a function $\mu$ from some Boolean algebra $\B$
into $[0,\infty]$ is a \emph{measure} if 
\begin{enumerate}[(a)]
  \item $\mu(0)=0$,
  \item\label{item:1} for every sequence $(a_n:n\in\N)$ of pairwise incompatible members of~$\B$
    (i.e.~$a_m\cdot a_n=0$ for $m\ne n$), if $\sum_{n=0}^\infty a_n\in\B$ 
    then $\mu\bigl(\sum_{n=0}^\infty a_n\bigr)=\sum_{n=0}^\infty\mu(a_n)$.
\end{enumerate}
Such a measure is \emph{strictly positive} if $a=0$ whenever $\mu(a)=0$, while a pair
$(\B,\mu)$ is a \emph{measure algebra} if $\B$ is $\sigma$-complete and $\mu$ is
strictly positive. A measure is \emph{totally finite} if $\mu(1)<\infty$. 

And $\mu:\B\to[0,\infty]$ is a \emph{finitely additive measure} if
\begin{enumerate}[(a)]
\item $\mu(0)=0$,
\item $\mu\bigl(\sum_{i=0}^{n-1}a_i\bigr)=\sum_{i=0}^{n-1}\mu(a_i)$ whenever 
  $a_0,\dots,a_{n-1}$ are pairwise incompatible.
\end{enumerate}

\subsection{Measuring over a model}
\label{sec:measuring-over-model}

To the author's knowledge, the following is a new notion. However, the special
case of a two-valued measure (i.e.~the range of $\mu$ is $\{0,1\}$) on an algebra of
subsets of $X$ is well-known for constructing elementary embeddings (e.g.~\cite{Lar}). 

\begin{defn}
Let $M$ be a class containing a set $X\in M$, and let $\A\subseteq M$ be a nonempty Boolean
algebra of subsets of $X$. A function $\mu:\A\to[0,\infty]$ is a
\emph{measure on $X$ over $M$} if 
\begin{enumerate}[(a)]
\item $\mu(\emptyset)=0$,
\item for every sequence $(E_n:n\in\N)$ of pairwise disjoint members of $\A$,
where $(E_n:n\in\N)\in M$, if $\bigcup_{n=0}^\infty E_n\in\A$ then 
$\mu\bigl(\bigcup_{n=0}^\infty E_n\bigr)=\sum_{n=0}^\infty\mu(E_n)$.
\end{enumerate}
And a \emph{measure space over $M$} is a triple $(X,\B,\mu)$ where $\B$ is a
$\sigma$-algebra of subsets of $X$ over $M$, i.e.~if $(E_n:n\in\N)\in M$ is a sequence of
members of $\B$ then $\bigcup_{n=0}^\infty E_n\in\B$, and $\mu:\B\to[0,\infty]$ is a
measure on $X$ over $M$.

As with the classical case, we generalize beyond algebras of subsets. Thus for a
Boolean algebra $\B\subseteq M$, $\mu:\B\to[0,\infty]$ is a \emph{measure over $M$} if
\begin{enumerate}[(a)]
\item $\mu(0)=0$,
\item for every sequence $(a_n:n\in\N)\in M$ forming an antichain, if
  $\sum_{n=0}^\infty a_n\in\B$ 
  then $\mu\bigl(\sum_{n=0}^\infty a_n\bigr)=\sum_{n=0}^\infty\mu(a_n)$. 
\end{enumerate}
\end{defn}

We will generally require the model $M$ to satisfy some small fragment of~$\zfc$,
and we write $\zfcf$ to denote enough as needed. In particular, 
we do not bother to generalize to a finitely additive measure over $M$, as we arrive
at the classical notion of finitely additive measures for these models.
We do not attempt any reverse
mathematical analysis, but do take note that in all of the results here 
neither the  power set axiom nor axiom of choice, beyond
countable choice, is necessary. 

Henceforth, we will focus exclusively on totally finite measures which form the most
important class. 

\begin{lem}\label{l-8}
Suppose $M$ is a model of $\zfcf$ containing a set $X\in M$, and $\mu:\A\to[0,\infty)$ is finitely
additive. If 
\textup{\begin{enumerate}[(A)]
\setcounter{enumi}{16-1}
\item\label{item:2} for every sequence $E_0\supseteq E_1\supseteq\cdots$ in $M$,  
$\inf_{n\in\N}\mu(E_n)=0$  whenever $\bigcap_{n=0}^\infty E_n=\emptyset$,
\end{enumerate}}
\noindent then $\mu$ is a measure over $M$.
\end{lem}
\begin{proof}
Let $(E_n:n\in\N)$ be a sequence in $M$ of pairwise disjoint members of $\A$, with
$E=\bigcup_{n=0}^\infty E_n\in\A$. First note that by finite additivity, 
\begin{equation}
  \mu(E)\ge\mu\left(\bigcup_{i=0}^{n-1}E_i\right)=\sum_{i=0}^{n-1}\mu(E_n)
\end{equation}
for all $n$. For the other inequality, suppose towards a contradiction that
$\mu(E)>\sum_{n=0}^\infty\mu(E_n)$. For each $n$, put $F_n=E\setminus\bigcup_{i=0}^{n-1}E_i$. Then
$F_0\supseteq F_1\supseteq\cdots$ and $\bigcap_{n=0}^\infty F_n=\emptyset$. Since moreover
the sequence $(F_n:n\in\N)$ is in $M$, $\inf_{n\in\N}\mu(F_n)=0$. And since
$\mu(E)=\mu(F_n)+\mu\bigl(\bigcup_{i=0}^{n-1}E_i\bigr)$ for all $n$, we have
\begin{equation}
  \mu(E)=\lim_{n\to\infty}\mu(F_n)+\sum_{i=0}^{n-1}\mu(E_i)=\sum_{n=0}^\infty \mu(E_n).\qedhere
\end{equation}
\end{proof}

We state without out proof the continuity property for measures over a model.

\begin{lem}
Let $M$ be a model of $\zfcf$, and let $\mu:\A\to[0,\infty)$ be a measure on $X$ over $M$. Then
\textup{\begin{enumerate}[(a)]
\item for every sequence $E_0\subseteq E_1\subseteq\cdots$ in $M$, 
of members of $\A$ where $\bigcup_{n=0}^\infty E_n\in\A$, 
$\mu\bigl(\bigcup_{n=0}^\infty E_n\bigr)=\sup_{n\in\N}\mu(E_n)$,
\item for every sequence $E_0\supseteq E_1\supseteq\cdots$ in $M$, of members of $\A$ where
  $\bigcap_{n=0}^\infty E_n\in\A$, $\mu\bigl(\bigcap_{n=0}^\infty E_n\bigr)=\inf_{n\in\N}\mu(E_n)$.
\end{enumerate}}
\end{lem}

The usual notion of an inner measure generalizes as follows. While outer measures
seem to be more popular, the (equivalent) use of inner measure merely reflects the
author's taste.

\begin{defn}
Let $M$ be a class containing a set $X\in M$. A function $\mu:\power(X)\cap M\to[0,\infty)$ is an
\emph{inner measure} on $X$ over $M$ if 
\begin{enumerate}[(a)]
\item $\mu(\emptyset)=0$,
\item $\mu(A\cup B)\ge\mu(A)+\mu(B)$ for all disjoint $A,B\in\power(X)\cap M$,
\item for every decreasing sequence $A_0\supseteq A_1\supseteq\cdots$ of members of
  $\power(X)\cap M$, where $(A_n:n\in\N)\in M$, 
  $\mu\bigl(\bigcap_{n=0}^\infty A_n\bigr)=\inf_{n\in\N}\mu(A_n)$.
\end{enumerate}
\end{defn}

Carath\'eodory's method for constructing a measure space from an inner measure on
some set does generalize to measuring over a model, as shown below. However, it is
of dubious value here, because nontrivial inner measures seem to be rather difficult
to construct. Thus while a given measure (in the classical sense) on an algebra of
sets can always be extended to a measure space, we  conjecture that, in general,
this is not true of measures over a model.

We recall this basic fact without giving proof.

\begin{lem}\label{l-2}
Suppose that $\A$ is an algebra of subsets of $X$ and $\rho:\A\to[0,\infty]$ is a
function with $\rho(\emptyset)=0$. Then
\begin{equation*}
  \E=\{E\in\A:\rho(A)=\rho(A\cap E)+\rho(A\setminus E)\textup{ for all }A\in\A\}
\end{equation*}
is an algebra of subsets of $X$, and $\rho$ is finitely additive on $\E$.
\end{lem}

\begin{thm}\label{T-1}
Let $\mu:\power(X)\cap M\to[0,\infty)$ be an inner measure over $M$, and set
\begin{multline*}
  \B=\{A\in\power(X)\cap M:
  \mu(E)=\mu(E\cap A)+\mu(E\setminus A)\\
  \textup{for all $E\in\power(X)\cap M$}\}.
\end{multline*}
Then $(X,\B,\mu\restriction\B)$ is a measure space over $M$.
\end{thm}
\begin{proof}
By Lemma~\ref{l-2}, $\B$ is an algebra of subsets of $X$ and $\mu\restriction\B$ is
finitely additive. Therefore, since $\mu$ is an inner measure over $M$, it follows from
Lemma~\ref{l-8} that $\mu\restriction\B$ is a measure over $M$.

Now suppose that $(A_n:n\in\N)\in M$ is a sequence of sets in $\B$. We need to show that
$A=\bigcup_{n=0}^\infty A_n\in\B$, so assume without loss of generality that the sequence is
increasing. Take $E\in\power(X)\cap M$. Since $E\setminus A_0\supseteq E\setminus
A_1\supseteq\cdots$, $\mu(E\setminus A)=\inf_{n\in\N}\mu(E\setminus A_n)$. 
Also note since $\mu$ is an inner measure,
$\mu(B)\le\mu(C)$ for all $B\subseteq C$ in $\power(X)\cap M$. Thus
\begin{equation}
  \mu(E)=\inf_{n\in\N}\mu(E\cap A_n)+\mu(E\setminus A_n)\le\mu(E\cap A)+\mu(E\setminus A).
\end{equation}
And the opposite inequality is immediate from the definition of an inner measure, proving
that $A\in\B$. 
\end{proof}

For any algebra $\A$ we let $\A_\delta$ denote the collection of
sets of the form $\bigcap_{n=0}^\infty E_n$ where $E_n\in\A$ for all $n$. Let
$\A_\sigma$ denote those sets whose complement is in $\A_\delta$, and
$\A_{\delta\cap\sigma}$ denote all sets of the form $A\cap B$ where $A\in\A_\delta$
and $B\in\A_\sigma$. 

Next we obtain a very slight extension for measures on some algebra $\A\subseteq M$ over 
$M$. Namely we extend, say $\mu:\A\to[0,\infty)$, to a finitely additive measure on
an algebra which includes $\A_\delta\cap M$, and which satisfies the natural
requirement
\begin{equation}
  \label{eq:26}
  \mu\left(\bigcap_{n=0}^\infty E_n\right)=\inf_{n\in\N}\mu(E_n)
\end{equation}
whenever $E_0\supseteq E_1\supseteq\cdots$ is a sequence in $M$ of members of $\A$. 
Notice that, even if we are only interested in finitely additive measures, by
Lemma~\ref{l-8}, we must begin with a measure over $M$ to satisfy this requirement.

The following lemma gives a standard representation of members of the algebra
generated by $\A_\delta$. The equivalent topological statement is: \emph{every set in
  the algebra generated by the open sets can be written as a finite disjoint union
  of sets of the form $F\cap U$ where $F$ is closed and $U$ is open}. 

\begin{lem}
\label{l-21}
Every set in the algebra generated by $\A_\delta$ is a finite disjoint union of
members of $\A_{\delta\cap\sigma}$. 
\end{lem}
\begin{proof}
First note that $\A_\delta$ is obviously closed under finite (or even countable)
intersections. And $\A_\delta$ is closed under finite unions because
$\bigcap_{n=0}^\infty E_n\cup\bigcap_{n=0}^\infty F_n
=\bigcap_{n=0}^\infty E_n\cup F_n$ whenever the sequences $(E_n:n\in\N)$ and
$(F_n:n\in\N)$ are decreasing, which we may assume for members of $\A_\delta$ as
$\A$ is an algebra.  Thus $\A_\sigma$ is also
closed under finite intersections and unions. 

To prove the lemma it suffices to show that the collection of finite disjoint unions
of members of $\A_{\delta\cap\sigma}$ is closed under finite intersections and
complements. But $\A_{\delta\cap\sigma}$ is closed under finite intersections by
prior considerations. And clearly the complement of an $\A_{\delta\cap\sigma}$ set
$A\cap B$, where $A\in\A_\delta$ and $B\in\A_\sigma$, 
can be written as the disjoint union of $A^\complement\cap B\in\A_\sigma$,
$A^\complement\cap B^\complement\in\A_{\delta\cap\sigma}$ and $A\cap
B^\complement\in\A_\delta$.
\end{proof}

Suppose $\mu:\A\to[0,\infty)$ is a measure over $M$.
Define $\bar\mu:\A_{\delta\cap\sigma}\cap M\to[0,1)$ by 
\begin{equation}
\bar\mu\left(\bigcap_{n=0}^\infty E_n\cap \bigcup_{n=0}^\infty F_n\right)
=\sup_{m\in\N}\inf_{n\in\N}\mu(E_n\cap F_m)    
\end{equation}
where $(E_n:n\in\N),(F_n:n\in\N)\in M$ are sequences of members of $\A$ such that
$E_0\supseteq E_1\supseteq\cdots$ and $F_0\subseteq F_1\subseteq\cdots$. Two
sequences with these properties are called a \emph{representation} of a member of
$\A_{\delta\cap\sigma}\cap M$.

\begin{notn}
  For a collection $\X$ of subsets of some base set $X$, let $\alg(\X)$ denote the
  algebra of sets generated by $\X$.
\end{notn}

We extend $\bar\mu$  to $\alg(\A_\delta)$ by letting 
\begin{equation}
\label{eq:16}
\bar\mu(E_0\cup\cdots\cup E_{n-1})=\sum_{i=0}^{n-1}\bar\mu(E_i)
\end{equation}
where the $E_i$'s are pairwise disjoint members of $\A_{\delta\cap\sigma}\cap M$. 

\begin{lem}
\label{l-10}
The mapping $\bar\mu:\alg(\A_\delta)\cap M\to[0,\infty)$ is a well-defined finitely
additive measure extending $\mu$.
\end{lem}
\begin{proof}
First of all we would like to be able to switch the order of $\sup$ and $\inf$.

\begin{claim}
\label{c-1}
Suppose that $(E_n,F_n:n\in\N)\in M$ is a representation of a member of $\A_{\sigma\cap\delta}$. 
Then $\inf_{m\in\N}\sup_{n\in\N}\mu(E_m\cap F_n)=0$
whenever $\bigcap_{n=0}^\infty E_n\cap\bigcup_{n=0}^\infty F_n=\emptyset$.
\end{claim}
\begin{proof}
Let $\varepsilon>0$ be arbitrary. Find $k$ so that 
\begin{equation}
  \label{eq:40}
  \mu(F_k)>\sup_{n\in\N}\mu(F_n)-\varepsilon.
\end{equation}
Then
\begin{equation}
  \label{eq:42}
  \begin{split}
  \inf_m\sup_n\mu(E_m\cap F_n)&=\inf_m\sup_n\mu(E_m\cap F_n\setminus F_k)
 +\mu(E_m\cap F_n\cap F_k)\\
 &\le\sup_n\mu(F_n\setminus F_k)+\inf_m\mu(E_m\cap F_k)\\
 &=\sup_n\mu(F_n)-\mu(F_k)+0\\
 &<\varepsilon.
 \end{split}
\end{equation}\\[-29pt]
\end{proof}

\begin{claim}
\label{c-3}
Suppose that $\bigcap_{n=0}^\infty E_n^i\cap\bigcup_{n=0}^\infty F_n^i$ \textup($i=0,\dots,j-1$\textup)
are pairwise disjoint members of $\A_{\delta\cap\sigma}\cap M$ whose union 
$\bigcap_{n=0}^\infty G_n\cap\bigcup_{n=0}^\infty H_n$ is a member of
$\A_{\delta\cap\sigma}$. Then 
\begin{equation*}
  \sum_{i<j}\inf_{m\in\N}\sup_{n\in\N}\mu(E_m^i\cap F_n^i)
  =\inf_{m\in\N}\sup_{n\in\N}\mu(G_m\cap H_n).
\end{equation*}
\end{claim}
\begin{proof}
Let us begin by noting that monotonicity implies 
\begin{equation}
  \label{eq:45}
  \inf_{m\in\N}\sup_{n\in\N}\mu(E^i_m\cap F^i_n)
  =\lim_{m\to\infty}\lim_{n\to\infty}\mu(E^i_m\cap F^i_n).
\end{equation}
For each $m,n$, for each $a\subseteq\{0,\dots,j-1\}$ 
set $J_{mn}^a=\bigcap_{i\in a}E_n^i\cap F_n^i$. It follows from Claim~\ref{c-1} and
pairwise disjointness, that for all $|a|\ge 2$, 
\begin{equation}
  \label{eq:38}
  \inf_m\sup_n \mu(J_{mn}^a)=0.
\end{equation}
Therefore, as
\begin{equation}
  \label{eq:32}
  \mu\Biggl(\bigcup_{i<j}E^i_m\cap F^i_n\Biggr)
  =\sum_{i<j}\mu(E^i_m\cap F^i_n)+\sum_{i=2}^{j-1}(-1)^{i-1}\cdot\sum_{a\in[j]^i}\mu(J_{mn}^a),
\end{equation}
we have
\begin{equation}
  \label{eq:43}
  \inf_m\sup_n\sum_{i<j}\mu(E^i_m\cap F^i_n)
  =\inf_m\sup_n\mu\Biggl(\bigcup_{i<j}E^i_m\cap F^i_n\Biggr).
\end{equation}
And for each $k,l,m$, the sequence $\bigcup_{i<j}E^i_k\cap F^i_l\setminus (G_m\cap H_0)
\supseteq\bigcup_{i<j}E^i_k\cap F^i_l\setminus (G_m\cap H_1)\supseteq\cdots$ is in $M$
and clearly has an empty intersection. Thus
\begin{equation}
  \label{eq:33}
  \inf_{n\in\N}\mu\Biggl(\bigcup_{i<j}E^i_k\cap F^i_l\setminus(G_m\cap H_n)\Biggr)=0.
\end{equation}

Now using~\eqref{eq:45},~\eqref{eq:43} and~\eqref{eq:33}, we obtain
\begin{equation}
  \label{eq:23}
  \begin{split}
  \sum_{i<j}\inf_{m}&\sup_n\mu(E^i_m\cap F^i_n)-\inf_m\sup_n\mu(G_m\cap H_n)\\
  &=\inf_m\sup_n\sum_{i<j}\mu(E^i_m\cap F^i_n)-\inf_m\sup_n\mu(G_m\cap H_n)\\
  &=\inf_m\sup_n\mu\Biggl(\bigcup_{i<j}E^i_m\cap F^i_n\Biggr)
    -\inf_m\sup_n\mu(G_m\cap H_n)\\
  &=\inf_{k}\sup_l\Biggl[\mu\Biggl(\bigcup_{i<j}E^i_k\cap F^i_l\Biggr)
  -\inf_m\sup_n\mu(G_m\cap H_n)\Biggr]\\
  &=\inf_k\sup_l\Biggl[\sup_m\Biggl(\mu\Biggl(\bigcup_{i<j}E^i_k\cap F^i_l\Biggr)
      -\sup_n\mu(G_m\cap H_n)\Biggr)\Biggr]\\
  &=\inf_k\sup_l\Biggl[\sup_m\inf_n\mu\Biggl(\bigcup_{i<j}E^i_k\cap F^i_l\Biggr)
      -\mu(G_m\cap H_n)\Biggr]\\
  &\le\inf_k\sup_l\Biggl[\sup_m\inf_n\mu\Biggl(\bigcup_{i<j}E^i_k\cap
     F^i_l\setminus(G_m\cap H_n)\Biggr)\Biggr]\\
  &=0
  \end{split}
\end{equation}
The other inequality is similar but easier. Using~\eqref{eq:45} and the appropriate
variation of~\eqref{eq:33}, we obtain
\begin{equation}
  \label{eq:21}
  \begin{split}
    \inf_m\sup_n\mu(G_m&\cap H_n)-\sum_{i<j}\mu(E^i_m\cap F_n^i)\\
    &=\inf_k\sup_l\Biggl[\sup_m\inf_n \mu(G_k\cap H_l)
    -\sum_{i<j}\mu(E_m^i\cap F_n^i)\Biggr]\\
    &\le\inf_k\sup_l\Biggl[\sup_m\inf_n\mu(G_k\cap H_l)
    -\mu\Biggl(\bigcup_{i<j}E_m^i\cap F_n^i\Biggr)\Biggr]\\
    &\le\inf_k\sup_l\Biggl[\sup_m\inf_n\mu\Biggr(G_k\cap H_l
    \setminus\bigcup_{i<j}E_m^i\cap F_n^i\Biggr)\Biggr]\\
    &=0.
  \end{split}
\end{equation}\\[-26pt]
\end{proof}

Claim~\ref{c-3} with $j=1$ states that $\bar\mu$ is well-defined on
$\A_{\delta\cap\sigma}\cap M$. And since any two `representations' 
of a set in $\alg(\A_\delta)$ as a finite disjoint union
of members of $\A_{\delta\cap\sigma}$ has a common refinement, it follows from
Claim~\ref{c-3} that $\bar\mu$ is well-defined. Finite additivity is an immediate
consequence of the definition.
\end{proof}

\subsection{Product measures}

The assumption that $\mu$ and $\nu$ are both measures over $M$ does not seem 
strong enough to guarantee the existence of a product measure over $M$. For this reason we consider
the notion of a measure over two models $M\subseteq N$ which is intermediate between 
a measure over $M$ and a measure over $N$.

\begin{defn}
For two models $M\subseteq N$, $X\in M$ and an algebra $\A\subseteq M$ of subsets of $X$,
a function $\mu:\A\to[0,\infty)$ is a \emph{measure over $M,N$} if
\begin{enumerate}[(a)]
\item $\mu(\emptyset)=0$,
\item for every sequence
$(E_n:n\in\N)\in M$ of pairwise disjoint members of $\A$, for every subsequence
$(E_{n_k}:k\in\N)\in N$ with $\bigcup_{k=0}^\infty E_{n_k}\in\A$,
$\mu\bigl(\bigcup_{k=0}^\infty E_{n_k}\bigr)=\sum_{k=0}^\infty\mu(E_{n_k})$.
\end{enumerate}
We say that $\mu$ is a \emph{measure allover $M$} if $\mu$ is a measure over $M,N$
for all $N\supseteq M$ satisfying $\zfcf$.
\end{defn}

\begin{lem}\label{l-15}
If $\mu:\A\to[0,\infty)$ is finitely additive, then $\mu$ is a measure over $M,N$ iff 
\textup{\begin{enumerate}[(A)]
\setcounter{enumi}{17-1}
\item\label{item:3} for every sequence $(E_n:n\in\N)\in M$, for every subsequence $E_{n_0}\supseteq
  E_{n_1}\supseteq\cdots$ in $N$, $\inf_{k\in\N}\mu(E_{n_k})=0$ whenever
  $\bigcap_{k=0}^\infty E_{n_k}=\emptyset$.
\end{enumerate}}
\end{lem}
\begin{proof}
Take a sequence $(E_n:n\in\N)\in M$ of pairwise disjoint members of $\A$, and suppose
$E_{n_0}\supseteq E_{n_1}\supseteq\cdots$ is a subsequence in $N$ with
$E=\bigcup_{k=0}^\infty E_{n_k}\in\A$. First note that by finite additivity, 
\begin{equation}
  \mu(E)=\mu\left(E\setminus\bigcup_{i=0}^{k-1}E_i\right)+\sum_{i=0}^{k-1}\mu(E_{n_i})
\end{equation}
for all $k$. Define for each $n$, $F_n=E\setminus \bigcup_{i=0}^{n-1}E_i$. Since the
sequence $(F_n:n\in\N)$ is in $M$, by property~\eqref{item:3}, $\inf_{k\in\N}\mu(F_{n_k})=0$. Thus
\begin{equation}
  \label{eq:25}
  \mu(E)=\lim_{k\to\infty}\mu(F_{n_k})+\sum_{i=0}^{k-1}\mu(E_{n_i})=\sum_{k=0}^\infty\mu(E_{n_k}).\qedhere
\end{equation}
\end{proof}

In what follows we take two algebras $\A$ and $\B$ of subsets of $X$ and $Y$,
respectively, and let $\A\otimes\B$ denote the algebra of subsets of $X\times Y$ consisting of all finite unions of sets of
the form $A\times B$ where $A\in\A$ and $B\in\B$. For a subset $G\subseteq X\times
Y$ and $x\in X$, we denote the projection of $G$ at $x$ by
\begin{equation}
  \label{eq:22}
  G_x=\{y\in Y:(x,y)\in G\}.
\end{equation}

If $G\in\A\otimes\B$, say $G=\bigcup_{i=0}^{n-1}E_i\times F_i$, then for all $x\in X$,
$G_x=\bigcup\{F_i:i<n$, $x\in E_i\}$, and therefore the set $\{G_x:x\in X\}$ has size at most
$2^n$ and in particular is finite. The \emph{standard representation} of $G$ is the
finite collection
\begin{equation}
  \label{eq:7}
  \stdrep(G)=\bigl\{\{x\in X:G_x=F\}\times F:F\in\{G_x:x\in X\}\bigr\}.
\end{equation}
Note that $G=\bigcup\stdrep(G)$.

For extended nonnegative real-valued functions $\mu$ on $\A$ and $\nu$ on $\B$,
define $\mu\otimes\nu:\A\otimes\B\to[0,\infty]$ by 
\begin{equation}
  \mu\otimes\nu(G)=\sum_{(E,F)\in\stdrep(G)}\mu(E)\nu(F)
\end{equation}
where $0\cdot\infty=0$. 

\begin{lem}\label{l-7}
Let $\mu:\A\to[0,\infty)$ and $\nu:\B\to[0,\infty)$ be functions, where $\mu$ is a finitely
additive measure on $X$. Then for all $G\in\A\otimes\B$ and all $0\le\sigma<\mu(X)$, 
\begin{equation*}
  \label{eq:35}
  \mu\bigl(\{x\in X:\nu(G_x)>\sigma\}\bigr)\ge\frac{\mu\otimes\nu(G)-\sigma}{\mu(X)-\sigma}.
\end{equation*}
\end{lem}
\begin{proof}
Put 
\begin{equation}
  \label{eq:8}
  E=\{x\in X:\mu(G_x)>\sigma\}.
\end{equation}
It is clear that $E\in\A$, and by the finite additivity of $\mu$,
$\mu\otimes\nu(G)\le(\mu(X)-\mu(E))\cdot\sigma+\mu(E)$. By rearranging terms, the result follows.
\end{proof}

Lemma~\ref{l-7} is often applied with $\sigma=1-\sqrt{1-\mu\otimes\nu(G)}$, in which
case the right hand side of~\eqref{eq:35} equals $\sigma$ when $\mu(X)=1$.

\begin{thm}\label{u-3}
Let $M\subseteq N$ be models of $\zfcf$ with $M$ transitive.
Suppose that $\A,\B\subseteq M$ are algebras of subsets of $X,Y\in M$, respectively.
Suppose further that $\mu:\A\to[0,\infty)$ is a measure over $M,N$, and that
$\nu:\B\to[0,\infty)$ is a measure over $M$ with $\nu\in N$.
Then $\mu\otimes\nu$ is a measure on $X\times Y$ over $M$.
\end{thm}
\begin{proof}
First we show that $\mu\otimes\nu$ is finitely additive. Note that if $\stdrep'(G)$ is a
finite refinement of $\stdrep(G)$ into pieces in $\A\times\B$, then
$\sum_{(E,F)\in\stdrep'(G)}\mu(E)\nu(F)=\mu\otimes\nu(G)$ by the finite additivity 
of $\mu$ and $\nu$. Suppose $G,H\in\A\otimes\B$ and
$G\cap H=\emptyset$. Then since $\stdrep(G)\cup\stdrep(H)$ is a refinement of
$\stdrep(G\cup H)$, we must have $\mu\otimes\nu(G\cup H)=\mu\otimes\nu(G)+\mu\otimes\nu(H)$.

Assume without loss of generality that $\mu(X)=1$.
Now by Lemma~\ref{l-8} it suffices to show that $\mu\otimes\nu$ satisfies property~\eqref{item:2}.
Suppose that $G_0\supseteq G_1\supseteq\cdots$ is a sequence in $M$ of members of
$\A\otimes\B$ where $\inf_{n\in\N}\mu\otimes\nu(G_n)=\delta>0$. By Lemma~\ref{l-7}, for
each $n$,
\begin{equation}
  \label{eq:29}
  \mu(E_n)\ge1-\sqrt{1-\delta},
\end{equation}
where $E_n=\bigl\{x\in X:\nu\bigl((G_n)_x\bigr)>1-\sqrt{1-\delta}\bigr\}$.
In $M$, choose a sequence $(A_n:n\in\N)$ which enumerates the set
\begin{equation}
  \label{eq:34}
  \bigcup_{n=0}^\infty\Bigl\{\bigcup\E:\E\subseteq\{E\in\A:(E,F)\in G_n\text{ for some $F\in\B$}\}\Bigr\}.
\end{equation}
Then, in $N$, $(E_n:n\in\N)$ can be written as a subsequence $(A_{n_k}:k\in\N)$ because
$\nu\in N$. Moreover, since the $E_n$'s are decreasing and since
$\inf_{k\in\N}\mu(A_{n_k})>0$, by the property~\eqref{item:3}, $\bigcap_{k=0}^\infty
A_{n_k}\ne\emptyset$. By the transitivity of $M$, $\bigcap_{k=0}^\infty
A_{n_k}\subseteq M$; hence, there is an $x\in X\cap M$ in the intersection. 
Then $(G_0)_x\supseteq(G_1)_x\supseteq\cdots$ is in $M$ and
$\inf_{n\in\N}\nu\bigl((G_n)_x\bigr)\ge1-\sqrt{1-\delta}>0$, and therefore there is a $y$ in
$\bigcap_{n=0}^\infty (G_n)_x$ because $\nu$ is a measure over $M$. 
Now $(x,y)$ witnesses that $\bigcap_{n=0}^\infty G_n\ne\emptyset$ as needed. 
\end{proof}

\section{The Measure}
\label{sec:themeas}

In this section we construct a measure over a certain model, which is of particular
interest to us.

\subsection{The forcing notion $\ctbl^+$}
\label{sec:forcing-notion-ctbl+}

For a set $S$, we let $\ctbl_S$ denote the ideal of all countable
subsets of $S$. Then $\ctbl_S^+$ is the coideal of all uncountable subsets of $S$; it
is considered as a forcing notion under the inclusion ordering.  We will primarily
be interested in $S=\oone$ in which case we just write $\ctbl^+$.

\begin{defn}
For a class $M$ and an index set $I\in M$, a filter $\F$ on $I$ is an \emph{ultrafilter over $M$} if
either $X\in\F$ or $I\setminus X\in\F$, for all $X\subseteq I$ in $M$. 
In the case $M=V$, we sometimes just say that $\F$ is an ultrafilter. 
A filter over $M$ is \emph{uniform} if all of its members have the same cardinality in $M$.
\end{defn}

For example, if $\U$ is a generic filter over $M$ of the poset $(\ctbl^+,\subseteq)$,
then $\U$ is generates a uniform ultrafilter on $\oone^M$ over $M$.

\begin{lem}
Let $\U$ be an ultrafilter on a set $S\in M$ over $M$. Then $f[[\U]]$ generates an
ultrafilter over $M$ for every $f:S\to S$ in $M$.
\end{lem}
\begin{proof}
Note that $f[[\U]]$ is a filter, because if $X,Y\in\U$ then $f[[\U]]\ni f(X\cap Y)\subseteq f(X)\cap f(Y)$.
Take $X\subseteq S$ in $M$. Then since $f\inv[X]\in M$, either
$f\inv[X]\in\U$ or $f\inv[S\setminus X]\in\U$.  Thus either
$f[f\inv[X]]\in f[[\U]]$ or $f[f\inv[S\setminus X]]\in
f[[\U]]$. We conclude that either $X\in \<f[[\U]]\>$ or $S\setminus
X\in \<f[[\U]]\>$. 
\end{proof}

\begin{lem}
\label{l-6}
Let $\U$ be $\ctbl^+$-generic over $M$.
If $f:\oone\to\oone$ is a countable--to--one function in $M$, then $f[[\U]]$ is
$\ctbl^+$ generic over $M$.
\end{lem}
\begin{proof}
Let $\D\subseteq\ctbl^+$ be a dense open subset in $M$. 
It suffices to show that the downwards closure of $f\inv[[\D]]$ is dense, 
because if $X\in\U\cap f\inv[[\D]]$,
then $f[X]\in f[[\U]]\cap \D$ since $f$ is countable--to--one and $\D$ is open.
And for all $X\in\ctbl^+$, there is a $Y\subseteq f[X]$ in $\D$, and then
$f\inv[Y]\cap X$ is in the downwards closure of $f\inv[[\D]]$.
\end{proof}

The following fact is well-known.

\begin{lem}\label{l-12}
Suppose that $\F$ is an ultrafilter on $I$ over $M$ and $f:I\to[-\infty,\infty]$ is a
function in $M$. Then the limit of $f$ at $\F$ exists, i.e.~there is an
$x\in[-\infty,\infty]$ such that $\lim_{i\to\F}f(i)=x$. Also for two extended
real-valued functions $f,g\in M$ on $I$, 
$\lim_{i\to\F}f(i)+g(i)=\lim_{i\to\F}f(i)+\lim_{i\to\F}g(i)$.
\end{lem}

\subsection{One dimensional measure}
\label{sec:one-dimens-meas}

A \emph{probability algebra} $(\random,\nu)$, i.e.~a measure algebra with
$\mu(1)=1$, is of interest to us as a forcing notion. Thus,
by Maharam's Theorem, we may assume that for some cardinal $\lambda$, 
$\random$ is the measure algebra of the measure space $2^\lambda$ where $\mu$ is the Haar
product measure. We view $\random$ as a defined notion. Thus, for example, if
$W\supseteq V$, we may have $W\models \random\ne\random^V$. For $a\in\random$ it
makes sense to define $\wh a$ as a coding of a member of $\random$ since every
member of $\random$ has a representative which is a Baire subset of $2^\lambda$. The
reader is warned that the `$\wh{\phantom a}$' is often suppressed. We write
$\C_\lambda\subseteq\random_\lambda$ for those elements which can be represented as
a clopen subset of $2^\lambda$. Recall that it forms a dense subset in the metric
topology on $\random_\lambda$: for all $a\in\random_\lambda$ and $\varepsilon>0$ there
exists $c\in\C_\lambda$ such that $\mu(a\diff c)<\varepsilon$.

The first Boolean algebra we consider is $(\random^{I};{\cup},{\cap},{\setminus})$,
where $I$ is some index set, 
with the operations inherited coordinatewise from $\random$, 
i.e.~$(\dot E\cup\dot F)(i)=\dot E(i)+\dot F(i)$, $(\dot E\cap \dot F)(i)=\dot E(i)\cdot \dot F(i)$ 
and $({\setminus}\dot E)(i)=-\dot E(i)$ for all $i\in I$. 
Members of $\random^I$ may be viewed as $\random$-names for subsets
of $I$, and of course the choice of operation symbols reflects the fact that 
\begin{equation}
  \label{eq:39}
  \lbrak\dot E\cup\dot F=(\dot E\cup\dot F)\rbrak
  =\lbrak\dot E\cap \dot F=(\dot E\cap\dot F)\rbrak
  =\lbrak\dot E\setminus\dot F=(\dot E\setminus\dot F)\rbrak
  =1
\end{equation}
where the set operation is intended on the left side of the equalities, and the
Boolean operation on the right. It is well-known that for an ultrafilter $\U$ on
$I$, $\lim_{i\to \U}(\dot E)$ defines a finitely additive measure on
$\random^I$. More familiarly, it defines a finitely additive measure (still not
necessarily strictly positive) on the ultrapower $\random^I\div\U$. The coding of
members of $\random$ naturally extends to members of $\random^I$, and we write
$\hat{\Dot E}$ to denote a coding of a member of $\random^I$, although we more often
suppress it.

For the purpose of stating the results, we fix two transitive models $M\subseteq N$
of $\zfcf$ and $\random=\random_\lambda$ for some $\lambda\in M$, 
but in usage, we will take $M=V$ and $N$ some forcing extension of $V$.
Let $I\in M$ be an index set, and let $\U$ generate an ultrafilter on~$I$
over $M$. In $N$, define
\begin{equation}
  \label{eq:47}
  \nu_\U:\random\times(\random^I)^M\to[0,1]
\end{equation}
by
\begin{equation}
  \label{eq:48}
  \nu_\U(a,\dot E)=\lim_{i\to\U}\mu\bigl(a\cdot\lbrak i\in\hat{\Dot E}\rmeas.
\end{equation}

\begin{lem}
The limit defining $\nu_\U$ always exists.
\end{lem}
\begin{proof}
Given $a\in\random$, let $(c_n:n\in\N)$ be a sequence in $\C_\lambda$ which converges
to $a$. Since each $c_n\in M$, by Lemma~\ref{l-12}, the limit defining $\nu_\U(\wh c_n,\dot E)$ exists. 
And since clearly $|\nu_\U(c,\dot E)-\nu_\U(d,\dot E)|\le\mu(c\diff d)$,
$(\nu_\U(c_n,\dot E):n\in\N)$ is a Cauchy sequence, and therefore has a convergent
subsequence, say indexed by $(n_k:k\in\N)$. 
It suffices to prove that
\begin{equation}
  \label{eq:24}
  \lim_{i\to\U}\mu\bigl(a\cdot\lbrak i\in\dot
  E\rbrak\bigr)=\lim_{k\to\infty}\nu_\U(c_{n_k},\dot E).
\end{equation}
Take $\varepsilon>0$.  Then choosing $k$ large enough, we have that
$\nu_\U(c_{n_k},\dot E)$ is within $\varepsilon\div3$ of the limit on the right and
$\mu(c_{n_k}\diff a)<\varepsilon\div 3$. Choose $X\in\U$ so that 
$\mu\bigl(c_{n_k}\cdot\lbrak i\in\dot E\rbrak\bigr)$ is within
$\varepsilon\div3$ of $\nu_\U(c_{n_k},\dot E)$ for all $i \in X$. 
Then $\mu\bigl(a\cdot\lbrak i\in\dot E\rbrak\bigr)$ 
is within $\varepsilon$ of the right hand limit for all $i\in X$, as needed.
\end{proof}

\begin{lem}\label{l-14}
The mapping $\nu_\U(\cdot,\dot E):\random\to[0,1]$ is a measure.
\end{lem}
\begin{proof}
Since it is a finite additive measure bounded by $\mu$.
\end{proof}

For $a\in\random^M$, the meaning of $(\random^I_a)^M$ is perfectly clear, but we would
rather consider $a\in\random$. Thus we associate an equivalence relation on
$(\random^I)^M$ with each $a\in\random$ in the natural manner: 
\begin{equation}
\label{eq:56}
\dot E\text{ is $a$-equivalent to $\dot F$}\Iff a\le\lbrak\dot E=\dot F\rbrak. 
\end{equation}
Then clearly the quotient map
\begin{equation}
  \label{eq:52}
  \bar\nu_\U(a,\cdot):(\random^I)^M\div a\to[0,1]
\end{equation}
(where $\bar\nu_\U(a,[\dot E]_a)=\nu_\U(a,\dot E)$) on the quotient Boolean algebra
$(\random^I)^M\div a$ is well-defined. 

For the remainder of this section, we are assuming that $M\models\ulc I$ is uncountable$\urc$.

\begin{lem}
\label{l-9}
Suppose $\U$ is $(\ctbl_I^+)$-generic over $M$. Then for all $a\in\random$, the
mapping $\nu_\U(a,\cdot):(\random^I)^M\div a\to[0,1]$ is a measure allover $M$. 
\end{lem}
\begin{proof}
It should be obvious that $\nu_\U(a,\emptyset)=0$. 
Finite additivity is immediate, 
because $\dot E\cap\dot F=0$ means that $\lbrak\dot E\cap\dot F=\emptyset\rbrak=1$.

Now we use Lemma~\ref{l-15}. Assume then that $W\supseteq M$ is a model of $\zfcf$,
$(\dot E_n:n\in\N)\in M$ is a sequence in $(\random^I)^M$, and there is a subsequence
$E_{n_0}\supseteq E_{n_1}\supseteq\cdots$ in $W$ for which $\bigcap_{k=0}^\infty E_{n_k}=\emptyset$. 

Take an arbitrary $\varepsilon>0$ in $M$. 
Choose $c\in\C_\lambda$ with $\mu(a\diff c)<\varepsilon\div2$. 
In $M$: We show that the set $\D$ of all conditions
of $\ctbl_I^+$ forcing that $\nu(\wh c,\dot E_{n_{\check k}})<\varepsilon\div 2$ for
some $k$ is dense. Take $X\in\ctbl_I^+$, and suppose towards a
contradiction that it has no extension in $\D$. 
Define $A\subseteq\N$ to consist
of all integers $n$ such that 
\begin{equation}
  \label{eq:27}
  Z_n=\bigl\{i\in X:\mu\bigl(c\cdot\lbrak i\in\dot E_n\rmeas<\varepsilon\div2\bigr\}
  \text{ is countable}.
\end{equation}

In $W$: By supposition, $\{n_k:k\in\N\}\subseteq A$. 
On the other hand, for each $i\in I$, 
$c\cdot\prod_{j=0}^\infty\sum_{k=j}^\infty\lbrak i\in\dot E_{n_k}\rbrak=0$, and
hence there is a $k_i$ such that 
\begin{equation}
\label{eq:49}
\mu\bigl(c\cdot\lbrak i\in\dot E_{n_{k_i}}\rmeas<\varepsilon.
\end{equation}
Therefore $\bigcup_{k=0}^\infty Z_{n_k}=X$. But this means that, in $M$,
$\bigcup_{n\in A}Z_n=X$, and we arrive at the absurdity that $M\models\ulc X$ is countable$\urc$.

Having proved that $\D$ is dense, by genericity, 
$\nu_\U(\wh c,\dot E_{n_k})<\varepsilon\div2$ for some~$k$, and thus $\nu_\U(a,\dot
E_{n_k})<\varepsilon$. This proves that the infimum is zero.
\end{proof}

\begin{defn}
If $\D$ is a directed partial ordering and $f:\D\to[-\infty,\infty]$, then for
$x\in[-\infty,\infty]$,
\begin{equation*}
  \lim_{p\to\D}f(p)=x
\end{equation*}
if
for every open $U\ni x$ there exists $p\in\D$ such that
$f[\D_p]\subseteq U$. Note that there is at most one such $x$.
\end{defn}

Continuing from~\eqref{eq:48}, suppose that $\H$ is $\random$-generic over
$N$. In $N[\H]$: Since $\random^M$ embeds into $\random$, there is a filter $\G\subseteq\random^M$
generic over $M$.   Define
\begin{equation}
m_\U:\power(I)^{M[\G]}\to[0,1]
\end{equation}
by
\begin{equation}
m_\U(E)=\lim_{a\to\H}\frac{\nu_\U(a,\dot E)}{\mu(a)},
\end{equation}
where $\dot E\in(\random^I)^M$ and $\dot E[\G]=E$. For the case of an ordinary
ultrafilter and $M=N$ and $\G=\H$ this is a familiar construction due to Solovay~\cite{Sol}.

\begin{lem}
\label{l-1}
\textup(In $N$\textup) for all $\delta\in[0,1]$, 
if $\{b\le a:\nu_\U(b,\dot E)\gtreqless\delta\cdot\mu(b)\}$ is dense below $a$, 
then $\nu_\U(a,\dot E)\gtreqless\delta$.
\end{lem}
\begin{proof}
Given $\rho<1$, find a finite antichain
$b_0,\dots,b_{k-1}$ below $a$ such that $\nu_\U(b_i,\dot E)\ge\delta\cdot\mu(b_i)$
for all $i=0,\dots,k-1$ and
\begin{equation}
  \sum_{i=0}^{k-1}\mu(b_i)>\sqrt\rho\cdot\mu(a).
\end{equation}
For each $i$ choose $X_i\in\U$ where $\mu\bigl(b_i\cdot\lbrak i\in\dot
E\rbrak\bigr)\ge\sqrt\rho\delta\cdot\mu(b_i)$ for all $i\in X_i$.
Put $X=\bigcap_{i=0}^{k-1}X_i\in\U$. Then for all $i\in X$,
\begin{equation}
\begin{split}
  \mu\bigl(a\cdot\lbrak i\in\dot E\rbrak\bigr)
&\ge\sum_{i=0}^{k-1}\mu\bigl(b_i\cdot\lbrak i\in\dot E\rbrak\bigr)
\ge\sum_{i=0}^{k-1}\sqrt{\rho}\delta\cdot\mu(b_i)\\
&\ge\rho\delta\cdot\mu(a).
\end{split}
\end{equation}
The other inequality is symmetrical.
\end{proof}

The following lemma is an immediate consequence.

\begin{lem}\label{l-4}
For all $\dot E\in(\random^I)^M$, $\delta\in[0,1]$ and $a\in\random^+$,
\begin{equation*} 
a\le\lbrak m_\U(\dot E)\gtreqless\delta\rbrak\Iff\nu_\U(b,\dot E)\gtreqless\delta\cdot\mu(b)
\textup{ for all }b\le a.
\end{equation*}
\end{lem}

\begin{lem}
$m_\U$ is a well-defined finitely additive measure on $I$.
\end{lem}
\begin{proof}
First we check that the limit indeed exists.
For each $\delta\in[0,1]$, put
\begin{equation}
  b_\delta=\sum\bigl\{a\in\random
    :\{b\le a:\nu_\U(b,\dot E)\ge\delta\cdot\mu(b)\}\text{ is dense below $a$}\bigr\}.
\end{equation}
Since $\{b\le b_\delta:\nu_\U(b,\dot E)\ge\delta\cdot\mu(b)\}$ is dense below
$b_\delta$, from Lemma~\ref{l-1} 
we conclude that $\nu_\U(c,\dot E)\ge\delta\cdot\mu(c)$ for all $c\le b_\delta$. Similarly,
$\nu_\U(c,\dot E)\le\delta\cdot\mu(c)$ for all $c\le-b_\delta$. 
Thus, by Lemma~\ref{l-4},  
\begin{equation}
  m_\U(\dot E)=\sup\{\delta:b_\delta\in\H\}.
\end{equation}

Now the fact that $m_\U$ is well-defined and finitely additive follows from Lemma~\ref{l-9}.
\end{proof}

\begin{lem}\label{l-13}
For all $0\le\sigma<1$,
\begin{equation*}
  \mu\bigl(a\cdot\lbrak m_\U(\dot E)>\sigma\rbrak\bigr)
  \ge\frac{\nu_\U(a,\dot E)-\sigma\cdot\mu(a)}{1-\sigma}.
\end{equation*}
\end{lem}
\begin{proof}
Put 
\begin{equation}
  \label{eq:41}
  b=a\cdot\lbrak m_\U(\dot E)\le\sigma\rbrak.
\end{equation}
Then by Lemma~\ref{l-4}, $\nu_\U(b,\dot E)\le\sigma\cdot\mu(b)$. Hence by the finite
additivity of $\nu_\U$ (Lemma~\ref{l-14}),
$\nu_\U(a,\dot E)=\nu_\U(a-b,\dot E)+\nu_\U(b,\dot E)\le\mu(a)-\mu(b)+\sigma\cdot\mu(b)$.
\end{proof}

\begin{lem}\label{l-3}
If $\U$ is $\ctbl_I^+$-generic over $M$, then $m_\U$ is a measure over $M[\G],N[\H]$.
\end{lem}
\begin{proof}
Since $m_\U$ is a finitely additive measure, by Lemma~\ref{l-15} it suffices to take a
sequence $(E_n:n\in\N)\in M[\G]$ of sets in
$\power(I)^{M[\G]}$, and a subsequence $E_{n_0}\supseteq E_{n_1}\supseteq\cdots$ in
$N[\H]$ where
$\bigcap_{k=0}^\infty E_{n_k}=\emptyset$, and prove that $\inf_{k\in\N}m_\U(E_{n_k})=0$. 
There is a sequence $(\dot E_n:n\in\N)\in M$ of members of $(\random^I)^M$ such that 
\begin{equation}
  \label{eq:58}
  \dot E_n[\G]=E_n\espc\text{for all $n$}.
\end{equation}
In $N[\H]$: There is an $a\in\H$ 
such that $a\le\lbrak\bigcap_{k=0}^\infty \dot E_{n_k}=\emptyset\rbrak$, or
equivalently $\prod_{k\in\N}[\dot E_{n_k}]_a=0$ in $(\random^I)^M\div a$. 
Hence by Lemma~\ref{l-9}, $\inf_{k\in\N}\nu_\U(a,\dot E_{n_k})=0$. This clearly implies
that $a\le \inf_{k\in\N}m_\U(E_{n_k})=0$, completing the proof.
\end{proof}

\begin{lem}
\label{l-17}
Let $f:\oone\to\oone$ be a function. Then for every $\dot E\in(\roone)^M$,
$m_{f[[\U]]}(\dot E)=m_\U(f\inv[\dot E])$.
\end{lem}
\begin{proof}
For all $a\in\random$,
\begin{equation}
\begin{split}
\lim_{\alpha\to f[[\U]]}\mu\bigl(a\cdot\lbrak\alpha\in\dot E\rbrak\bigr)
&=\lim_{\alpha\to \U}\mu\bigl(a\cdot\lbrak f(\alpha)\in\dot E\rbrak\bigr)\\
&=\lim_{\alpha\to\U}\mu\bigl(a\cdot\lbrak\alpha\in f\inv[\dot E]\rbrak\bigr),
\end{split}
\end{equation}
and thus $\nu_{f[[\U]]}(a,\dot E)=\nu_{\U}(a,f\inv[\dot E])$.
\end{proof}

\subsection{Product measure}
\label{sec:product-measure}

Now we arrive at the construction of the primary and secondary measures ($m$ and
$\nu$, resp.). Let $\U$ be $\ctbl_I^+$-generic over $M$, let $\V$ be $(\ctbl_I^+)^M$-generic over
$M[\U]$ and let $\H$ be $\random^{M[\U][\V]}$-generic over $M[\U][\V]$.
In other words $M[\U][\V][\H]$ is a generic extension of $V$ by the poset
\begin{equation}
\label{eq:50}
(\ctbl_I^+\times\ctbl_I^+)\gdot\random.
\end{equation}

By Lemma~\ref{l-9}, for all $a\in\random$, the restriction
$\nu_\U(a,\cdot):(\random^I)^M\div a\to[0,1]$ is a measure over $M,M[\U][\V]$; and
also, $\nu_\V(a,\cdot):(\random^I)^M\div a\to[0,1]$ is a measure over $M$.
Therefore, by Theorem~\ref{u-3}, 
\begin{equation}
  \label{eq:51}
  \nu_\U(a,\cdot)\otimes\nu_\V(b,\cdot):(\random^I)^M\div a\otimes(\random^I)^M\div b\to[0,1]
\end{equation}
is a measure over $M$, for each $a,b\in\random$. And by Lemma~\ref{l-10}, 
it has an extension 
\begin{equation}
  \label{eq:63}
  \overline{\nu_\U(a,\cdot)\otimes\nu_\V(b,\cdot)}:
  \alg\bigl((\random^I)^M\div a\otimes (\random^I)^M\div b\bigr)_\delta\cap M\to[0,1]
\end{equation}
to a finitely additive measure such that 
\begin{equation}
  \label{eq:77}
  \overline{\nu_\U(a,\cdot)\otimes\nu_\V(b,\cdot)}\left(\bigcap_{n=0}^\infty G_n\right)
  =\inf_{n\in\N}\nu_\U(a,\cdot)\otimes\nu_\V(b,\cdot)\bigl(G_n\bigr)
\end{equation}
for every sequence $G_0\supseteq G_1\supseteq\cdots$ in $M$ of members 
of $(\random^I)^M\div a\otimes (\random^I)^M\div b$. In the case  $a=b$, we denote
the function in~\eqref{eq:63} by $\nu(a,\cdot)$. 

By Lemma~\ref{l-3}, $m_\U:\power(I)^{M[\G]}\to[0,1]$ defined by
\begin{equation}
  \label{eq:1}
  m_\U(E)=\lim_{a\to\H}\lim_{i\to\U}\frac{\mu(a\cdot\| i\in\dot E\|)}{\mu(a)}
  \espc\text{where $\dot E[\G]=E$}
\end{equation}
is a measure over
$M[\G],M[\U][\V][\H]$; and $m_\V:\power(I)^{M[\G]}\to[0,1]$ is a measure over
$M[\G]$. Thus $m_\U\otimes m_\V:\power(I)^{M[\G]}\otimes\power(I)^{M[\G]}\to[0,1]$
is a measure over $M[\G]$, and it has an extension, which we denote by $m$, to
\begin{equation}
  \label{eq:17}
  m:\alg\bigl(\power(I)^{M[\G]}\otimes\power(I)^{M[\G]}\bigr)_\delta\cap M\to[0,1]
\end{equation}
satisfying the continuity condition~\eqref{eq:26}. 

Henceforth, we will be working in $V$, except where otherwise noted, and thus the
$M$ symbols above vanish.

The following Lemma should be able to be improved to a more informative
estimate. Nonetheless, even this is nontrivial to prove, because it requires
Fubini's Theorem to be generalized to our context; and, due to time constraints, this
will not be done in the present note.

\begin{lem}
Let $\dot K\in \alg(\roone_a\otimes\roone_a)_\delta$ for some $a\in\random$. For all
$(X,Y)\in\ctbl^+$, if
\begin{equation*}
  (X,Y)\forces\wh a \le\lbrak\nu(\dot K)=1\rbrak,
\end{equation*}
then for all $\rho<1$ there exists $\alpha\in X$ and $\beta>\alpha$ in $Y$ such that 
\begin{equation*}
  \mu\bigl(a\cdot\lbrak\dot K(\alpha,\beta)\rmeas>\rho\cdot\mu(a).
\end{equation*}
\end{lem}
\begin{proof}
Omitted.
\end{proof}

\subsection{Separable measure algebra}
\label{sec:separ-meas-algebra}

\begin{thm}
If $\U$ is $\ctbl^+$ generic over $V$, and $\random$ is a separable measure algebra,
then $m_\U$ is a two-valued measure, i.e.~it has range $\two$.
\end{thm}
\begin{proof}
Let $\dot E\in V$ be an $\random$-name for a member of $\power(\oone)$. In $V[\U]$:
Supposing $a\in\random^+$ forces that $m_\U(\dot E)>0$, we find an extension of $a$ forcing
$m_\U(\dot E)=1$. Now there is a $b\le a$ and a rational $\varepsilon>0$ such that 
\begin{equation}
  \label{eq:5}
  b\le\lbrak m_\U(\dot E)\rbrak\ge\varepsilon,
\end{equation}
and there is a rational $\sigma>0$ such that $\mu(b)\ge\sigma$.
Hence by Lemma~\ref{l-4}, $\nu_\U(b,\dot E)\ge\varepsilon$ and thus there is an
$X\in\U$ such that 
\begin{equation}
  \label{eq:65}
  \mu\bigl(b\cdot\lbrak\alpha\in\dot E\rmeas>\frac{3\varepsilon\sigma}4\espc\text{for all
    $\alpha\in X$}.
\end{equation}
Take $c\in\C$ such that 
\begin{equation}
  \label{eq:67}
  \mu(b\diff c)<\frac{(1-\rho)\varepsilon\sigma}4.
\end{equation}
Then in particular, 
\begin{equation}
  \label{eq:66}
  \mu\bigl(c\cdot\lbrak\alpha\in\dot E\rmeas
  >\frac{\varepsilon\sigma}2\espc\text{for all $\alpha\in X$}.
\end{equation}
Therefore, for each $\alpha\in X$, there exists $d_\alpha\le c$ in $\C$ such that
\begin{align}
  \label{eq:68}
  \mu(d_\alpha)&\ge\frac{\varepsilon\sigma}2,\\
  \label{eq:69}
  \mu\bigl(d_\alpha\diff\bigl(c\cdot\lbrak\alpha\in\dot E\rmeas\bigr)&<\frac{(1-\rho)\varepsilon\sigma}4.
\end{align}
A computation with~\eqref{eq:67}, \eqref{eq:68} and~\eqref{eq:69} yields
\begin{equation}
  \label{eq:70}
  \mu\bigl(b\cdot d_\alpha\cdot\lbrak\alpha\in \dot E\rbrak\bigr)
  >\mu(d_\alpha)-\frac{(1-\rho)\varepsilon\sigma}2\ge\rho\cdot\mu(b\cdot d_\alpha).
\end{equation}
Since $(d_\alpha:\alpha\in X)\in V$, by genericity, there exists $X_0\subseteq X$
in $\U$ such that $d_\alpha$ equals some fixed $d$ for all $\alpha\in X_0$. This
proves that $\nu_\U(b\cdot d,\dot E)\ge\rho\cdot\mu(b\cdot d)$, as required.
\end{proof}

\begin{lem}
\label{l-5}
Let $\dot K\in\alg\bigl((\otimes_{i=0}^{n-1}\roone)_\delta\bigr)$ for some
$n=1,2,\dots$. Then 
\begin{equation}
  \label{eq:10}
  \bigcap_{\varepsilon>0}\{a\in\random:\vector X\forces\nu(\wh a,\dot K)<\varepsilon
  \textup{ for some }\vector X\in(\ctbl^+)^n\}
\end{equation}
is dense below
\begin{equation}
  \label{eq:12}
  -\sum\bigl\{a\in\random:(\ctbl^+)^n\forces\wh a\le\lbrak m(\dot K)=1\rbrak\bigr\}.
\end{equation}
\end{lem}
\begin{proof}
Let $\dot a$ be a $\ctbl^+$-name for $\lbrak m(\dot K)=0\rbrak$. 
Take a nonzero $c$ beneath the Boolean value in~\eqref{eq:12}. 
Then for some $\vector X\in\ctbl^+$,
\begin{equation}
  \label{eq:60}
  \vector X\forces\dot a\cdot\wh c\ne0.
\end{equation}
Find $\vector X_0\le \vector X$ and $\delta>0$ where
\begin{equation}
  \label{eq:3}
  \vector X_0\forces \mu(\dot a\cdot c)>\delta.
\end{equation}
Choose a sequence $\vector X_0\ge \vector X_1\ge\cdots$ of conditions and $d_i\in\random$ such that 
\begin{equation}
  \label{eq:6}
  \vector X_i\forces \mu\bigl((\dot a\cdot c)\diff d_i\bigr)<\frac{\delta}{2^i}
  \espc\text{for all $i=1,2,\dots$}.
\end{equation}
Put $d=\sum_{k=0}^\infty\prod_{i=k}^\infty c\cdot d_i$. 
Since $\vector X_i\forces\mu(c\cdot d_i)>\mu(\dot a\cdot c)-\delta\div2^i$, by~\eqref{eq:3},
\begin{equation}
  \label{eq:61}
  \mu(c\cdot d_i)>(1-2^{-i})\delta\espc\text{for all $i$,}
\end{equation}
and thus $\mu(d)\ge\delta$. 
And since $\ctbl^+\forces\nu(\dot a,\dot K)=0$, 
$\vector X_i\forces\nu(\wh d,\dot E)<\delta\div 2^i$, completing the proof that $d$ is
a nonzero member of the set~\eqref{eq:10} below~$c$.
\end{proof}

\subsection{Measure homogeneity}
\label{sec:measure-homogeneity}

\begin{defn}
\label{d-1}
Given a subfamily $\F\subseteq[\oone]^n$ for some $n=1,2,\dots$, let
$f_i:\oone\to\oone$ ($i=0,\dots,n-1$) be the nondecreasing functions
where $\F=\bigl\{\{f_0(\alpha)<\dots<f_{n-1}(\alpha)\}:\alpha<\oone\bigr\}$. 
Then element $\dot K\in\random^{\oone^2}$ of the random power determines 
a colouring $\dot K_\F\in\random^{\oone^2}$ via
\begin{multline}
  \label{eq:9}
  \lbrak \dot K_\F(\alpha,\beta)\rbrak
  =\prod_{i=0}^{n-1}\prod_{j=0}^{n-1}-\lbrak\dot K\bigl(f_i(\alpha),f_j(\alpha)\bigr)\rbrak\\
  \cdot\prod_{i=0}^{n-1}\prod_{j=0}^{n-1}-\lbrak\dot K\bigl(f_i(\beta),f_j(\beta)\bigr)\rbrak\\
  \cdot\sum_{i=0}^{n-1}\sum_{j=0}^{n-1}\lbrak K\bigl(f_i(\alpha),f_j(\beta)\bigr)\rbrak.
\end{multline}
\end{defn}

The following notion, while not particularly natural, will serve for the purposes of
this note.

\begin{defn}
An $\alg(\roone\otimes\roone)_\delta$ set $\dot K\in\random^{\oone}$
random colouring $\dot K\in\random^{\oone^2}$ is called
\emph{$\ctbl^+\times\ctbl^+$-measure zero homogeneous at $a\in\random$} 
if there exists a dense set of elements $b\le a$ such that: 
for every $n=1,2,\dots$, for every uncountable pairwise
disjoint subfamily $\F\subseteq[\oone]^n$, for all $\varepsilon>0$, there exists
$(X,Y)\in\ctbl^+\times\ctbl^+$ forcing that
\begin{equation*}
  \nu(\wh b,\dot K_\F)<\varepsilon.
\end{equation*}
\end{defn}

\section{Forcing with a measure algebra}

Recall that for a given subset $K\subseteq S^2$, 
a subset $X\subseteq S$ is called \emph{$K$-homogeneous} if $(\alpha,\beta)\in K$ for
all $\alpha<\beta$ in $X$. 

The following Theorem is the main result for applying the measure constructed in
Section~\ref{sec:themeas} to the analysis of measure algebraic forcing.

\begin{thm}[$\m>\aleph_1$]\label{u-2}
Let $\random$ be any measure algebra. Suppose that
$\dot K$ is an $\random$-name for a subset of $\oone^2$ in
$\alg\bigl(\power(\oone)\otimes\power(\oone)\bigr)_\delta$, and $a\in\random$.
If 
\begin{equation*}
\ctbl^+\times\ctbl^+\forces \wh a\le\lbrak m(\dot K)=1\rbrak,
\end{equation*} 
then $a$ forces the existence of a countable decomposition of $\oone$
into $\dot K$-homogeneous pieces.
\end{thm}

\begin{lem}\label{l-11}
Let $\dot K$ be an $\random$-name for a member of 
$\alg\bigl(\power(\oone)\otimes\power(\oone)\bigr)_\delta$ and $a\in\random$, and suppose
$(\ctbl^+\times\ctbl^+)\forces\wh a\le\lbrak m(\dot K)=1\rbrak$.
If $\F$ is an uncountable pairwise disjoint family of finite subsets of $\oone$, 
then for every $\rho<1$ there are
$\vgamma<\vlambda$ in $\F$ such that
\begin{equation*}
  \mu\Biggl(a\cdot\prod_{\alpha\in\vgamma}\prod_{\beta\in\vlambda}\lbrak\dot K(\alpha,\beta)\rbrak\Biggr)
  >\rho\cdot\mu(a).
\end{equation*}
\end{lem}
\begin{proof}
Without loss of generality assume $a=1$.
Given $\F$ we may as well assume that for some $n\ge 1$, $|\vgamma|=n$ for all
$\vgamma\in\F$. Let $f_i:\oone\to\oone$ for $i=0,\dots,n-1$ be functions where
$\F=\bigl\{\{f_0(\alpha),\dots,f_{n-1}(\alpha)\}:\alpha<\oone\bigr\}$. 
For each $i,j=0,\dots,n-1$, define an $\random$-name $\dot K_{ij}$ for
  a relation on $\oone$ by
  \begin{equation}
    \lbrak\dot K_{ij}(\alpha,\beta)\rbrak=\lbrak\dot K\bigl(f_i(\alpha), f_j(\beta)\bigr)\rbrak.
  \end{equation}
As follows from Lemmas~\ref{l-6} and~\ref{l-17}, 
$\pob$ forces $m\bigl(\bigcap_{i,j=0}^{n-1}\dot K_{ij}\bigr)=1$. 
And Lemma~\ref{l-5} gives $\alpha<\beta$ 
such that $\lmeas(\alpha,\beta)\in\bigcap_{i,j=0}^{n-1}\dot K_{ij}\rmeas>\rho$ as required.
\end{proof}

\begin{proof}[Proof of Theorem~\textup{\ref{u-2}}]
For $\varepsilon>0$, let $\Q^1_\varepsilon$ be the poset of all finite $\vgamma\subseteq\oone$ such that 
\begin{equation}
  \label{eq:36}
  \mu\left(\sum_{\alpha\in\vgamma\cap\beta}
   -\lbrak\dot K(\alpha,\beta)\rbrak\right)<\varepsilon\espc\text{for all $\beta\in\vgamma$},
\end{equation}
ordered by set containment.

\begin{claim}
$\Q^1_\varepsilon$ has the ccc.
\end{claim}
\begin{proof}
Supposing $\Scal_0$ is an uncountable subset of $\Q^1_\varepsilon$, 
by going to an uncountable subset we
can assume that $\Scal_0$ is a $\Delta$-system, say with root $\vlambda$, and that for some
integer $n$, $|\vgamma\setminus\vlambda|=n$ for every  $\vgamma\in\Scal_0$. Find an
uncountable $\Scal_1\subseteq\Scal_0$ and $\tau<\varepsilon$ such that
\begin{equation}
  \label{eq:37}
  \mu\left(\sum_{\alpha\in\vgamma\cap\beta}-\lbrak\dot
  K(\alpha,\beta)\rbrak\right)<\tau\espc\text{for all $\beta\in\vgamma\setminus\vlambda$},
\end{equation}
for all $\vgamma\in\Scal_1$. By Lemma~\ref{l-11}, there are $\vgamma_0,\vgamma_1$ in $\Scal_1$
such that $(\vgamma_0\setminus\vlambda)<(\vgamma_1\setminus\vlambda)$ and 
$\prod_{\alpha\in\vgamma_0\setminus\vlambda}\prod_{\beta\in\vgamma_1\setminus\vlambda}
\lbrak\dot K(\alpha,\beta)\rbrak$ has measure greater than $1-(\varepsilon-\tau)$. This implies in
particular that
\begin{equation}
  \label{eq:44}
  \mu\left(\sum_{\alpha\in\vgamma_0\setminus\vlambda}
  -\lbrak\dot K(\alpha,\beta)\rbrak\right)<\varepsilon-\tau
  \espc\text{for all $\beta\in\vgamma_1\setminus\vlambda$}.
\end{equation}
Hence for every $\beta\in\vgamma_1\setminus\vlambda$,
\begin{equation}
\mu\left(\sum_{\alpha\in(\vgamma_0\cup\vgamma_1)\cap\beta}-\lbrak\dot K(\alpha,\beta)\rbrak\right)
\le\tau+\mu\left(\sum_{\alpha\in\vgamma_0\setminus\vlambda}
-\lbrak\dot K(\alpha,\beta)\rbrak\right)<\varepsilon, 
\end{equation}
proving that $\vgamma_0\cup\vgamma_1\in\Q^1_\varepsilon$.
\end{proof}

Suppose that $\C\subseteq\Q^1_\varepsilon$ is centered, and put $X=\bigcup\C$. Then define
an $\random$-name $\dot Y$ for a subset of $\oone$ by
\begin{equation}
  \label{eq:46}
  \lbrak\beta\in\dot Y\rbrak=\prod_{\alpha\in X\cap\beta}\lbrak\dot K(\alpha,\beta)\rbrak
\end{equation}
for all $\beta\in X$. It is clear that $\dot Y$ is $\dot K$-homogeneous with probability
one, and that $\lmeas\beta\in\dot Y\rmeas\ge1-\varepsilon$ for all $\beta\in X$. Thus a
straightforward density argument applied to the finite support iteration of the sequence
of posets $\bigl(\Q^1_{(n+1)\inv}:n\in\N\bigr)$ gives the desired $\random$-name for a countable
decomposition of $\oone$.
\end{proof}

\begin{thm}[$\m>\aleph_1$]
\label{u-4}
Suppose $\random$ is a separable measure algebra.
If $\dot K\in(\roone\otimes\roone)_\delta$ is measure zero
homogeneous at $a\in\random$, then $a$ forces that $\oone$ has a countable
decomposition into $\oone^2\setminus\dot K$-homogeneous pieces.
\end{thm}
\begin{proof}
We may assume $a=1$ without loss of generality.
The poset $\Q^0_\varepsilon$ consists of all finite  $\varGamma\subseteq\oone$ such that 
\begin{equation}
  \label{eq:11}
  \mu\left(\sum_{\substack{ \alpha,\beta\in\varGamma\\ \alpha<\beta}}
    \lbrak \dot K(\alpha,\beta)\rbrak\right)<\varepsilon.
\end{equation}

\begin{claim}
$\Q^0_\varepsilon$ has the ccc.
\end{claim}
\begin{proof}
Given an uncountable subsets $\Scal\subseteq\Q^0_\varepsilon$ we can assume that $\Scal$
is a $\Delta$-system, say with root $\varLambda$, such that
$|\varGamma\setminus\varLambda|=n$ for all $\varGamma\in\Scal$. By further refinement,
assume that 
\begin{equation}
  \label{eq:13}
  \mu\left(\sum_{\substack{\alpha,\beta\in\varGamma\\ \alpha<\beta}}\lbrak\dot
    K(\alpha,\beta)\rbrak\right)<\tau_0
  \espc\text{for all $\varGamma\in\Scal$}
\end{equation}
for some $\tau_0<\varepsilon$. Then by separability there exists an uncountable
$\Scal_0\subseteq\Scal$ such that 
\begin{equation}
  \label{eq:14}
  \mu\left(\sum_{\substack{\alpha,\beta\in\varGamma_0\\ \alpha<\beta}}\lbrak\dot K(\alpha,\beta)\rbrak
    +\sum_{\substack{\alpha,\beta\in\varGamma_1\\ \alpha<\beta}}\lbrak\dot K(\alpha,\beta)\rbrak\right)
  <\tau_1
  \espc\text{for all $\varGamma_0,\varGamma_1\in\Scal_0$}
\end{equation}
for some $\tau_1<\varepsilon$. Therefore, it suffices to find
$\varGamma_0<\varGamma_1$ in $\Scal_0$ such that 
\begin{multline}
  \label{eq:15}
  \mu\left(\vphantom{\sum_{\beta\in\varGamma}}\right.
  \prod_{\substack{\alpha,\beta\in\varGamma_0\\ \alpha<\beta}}-\lbrak\dot K(\alpha,\beta)\rbrak
    \cdot\prod_{\substack{\alpha,\beta\in\varGamma_1\\ \alpha<\beta}}-\lbrak\dot K(\alpha,\beta)\rbrak
   \\ \left.
    \cdot\sum_{\alpha\in\varGamma_0\setminus\varLambda}
    \sum_{\beta\in\varGamma_1\setminus\varLambda}\lbrak\dot K(\alpha,\beta)\rbrak\right)
  <\varepsilon-\tau_1.
\end{multline}
And their existence follows from the hypothesis that some $(X,Y)$ forces
$\nu(1,\dot K_\F)<\varepsilon-\tau_1$ where
$\F=\{\varGamma\setminus\varLambda:\varGamma\in\Scal_0\}$. 
\end{proof}

Now precisely the same argument as in the proof of Theorem~\ref{u-2} yields the
desired countable decomposition.
\end{proof}

\section{The measure of an $(\oone,\oone)$ gap.}
\label{sec:measure-an-oone}

The reader is referred to~\cite{Hir} for the notion of a (Souslin) pregap in $L^0(\nu)\div\Fin$.
Let $(\dot x_\alpha,\dot y_\alpha:\alpha<\oone)$ be sequences of
representatives of an $(\oone,\oone)$-pregap in $L^0(\nu)\div\Fin$; in other words,
the $\dot x_\alpha$'s and $\dot y_\alpha$'s are $\random$-names for subsets of $\N$,
such that the pair of sequences name a pregap in $\pN\div\Fin$. 
Then for each $n\in\N$, $\dot E_n\in\roone$ and $\dot F_n\in\roone$ are determined by
\begin{equation}
  \label{eq:75}
  \dot E_n=\{\alpha<\oone:n\in\dot x_\alpha\}\And\dot F_n=\{\beta<\oone:n\in\dot y_\beta\}.
\end{equation}
This associates with the pregap an $\random$-name for a subset of $\oone^2$
\begin{equation}
\label{eq:28}
\dot K\in\bigl(\power(\oone)\otimes\power(\oone)\bigr)_\sigma
\end{equation}
via
\begin{equation}
  \label{eq:76}
  \dot K=\bigcup_{n=0}^\infty \dot E_n\times\dot F_n.
\end{equation}
Thus
\begin{equation*}
  \lbrak\dot K(\alpha,\beta)\rbrak=\lbrak\dot x_\alpha\cap\dot y_\beta\ne\emptyset\rbrak
  \espc\text{for all $\alpha,\beta$}.
\end{equation*}

The random colouring $\dot K$ in turn determines a member of $\random$ by
\begin{equation}
  \label{eq:31}
  b=\sum\bigl\{a\in\random:\ctbl^+\times\ctbl^+\forces\wh a\le\lbrak m(\dot K)=1\rbrak\bigr\}.
\end{equation}

\begin{lem}
\label{l-16}
If $\random$ is separable, then for every $(\oone,\oone)$ pregap the associated
colouring $\dot K$ is measure zero homogeneous at $-b$ as specified in~\eqref{eq:31}.
\end{lem}
\begin{proof}
Let $\F\subseteq[\oone]^n$ be a given uncountable pairwise disjoint family, and let
$f_i:\oone\to\oone$ ($i=0,\dots,n-1$) be the nondecreasing functions where
$\F=\bigl\{\{f_0(\alpha)<\cdots<f_{n-1}(\alpha)\}:\alpha<\oone\bigr\}$.
Take $a\le -b$. By Lemma~\ref{l-5} there exists $c\le a$ such that 
\begin{equation}
  \label{eq:53}
  (X,Y)\forces\nu(\wh c,\dot K)<\varepsilon\espc\text{for some $(X,Y)\in(\ctbl^+)^2$}
\end{equation}
for all $\varepsilon>0$. Now fix $\varepsilon>0$. Then choose $(X,Y)$ forcing that
\begin{equation}
  \label{eq:55}
  \nu(\wh c,\dot K)<\frac\varepsilon2.
\end{equation}
Choose functions $g:\oone\to X$ and $h:\oone\to Y$ such that 
\begin{alignat}{2}
  \label{eq:71}
  g(\alpha)&\ge f_{n-1}(\alpha)\espc&&\text{for all $\alpha<\oone$},\\
  \label{eq:72}
  h(\beta)&\ge f_{n-1}(\beta)&&\text{for all $\beta<\oone$}.
\end{alignat}

\begin{claim}
\label{c-4}
$(Z,Z)\forces\nu\bigl(\wh c,\dot K\circ(f,g)\bigr)<\varepsilon\div2$, 
i.e.~$\dot K\circ(f,g)(\alpha,\beta)=\dot K(f(\alpha),g(\beta))$.
\end{claim}
\begin{proof}
Let $(\U,\V)$ be $\ctbl^+\times\ctbl^+$-generic with $Z\in\U,\V$. Then by
Lemma~\ref{l-6}, $(f[[\U]],g[[\V]])$ is $\ctbl^+\times\ctbl^+$-generic. In $V[\U][\V]$:
And since $(f[Z],g[Z])\le(X,Y)$, by~\eqref{eq:55} and Lemma~\ref{l-17},
\begin{equation}
  \begin{split}
  \label{eq:18}
  \nu\bigl(c,(E_n\times F_n)\circ(f,g)\bigr)
  &=\nu_\U(c,f\inv[\dot E_n])\nu_\V(c,g\inv[\dot F_n])\\
  &=\nu_{f[[\U]]}(c,\dot E_n)\nu_{g[[\V]]}(c,\dot F_n)\\
  &<\frac\varepsilon2
  \end{split}
\end{equation}
for all $n$. Since $\nu(c,\dot K\circ(f,g))
=\sup_{n\in\N}\nu\bigl(c,(E_n\times F_n)\circ(f,g)\bigr)$ the proof is complete.
\end{proof}

Choose $k$ large enough so that
\begin{alignat}{2}
  \label{eq:57}
  \lbrak\dot x_{f_i(\alpha)}\setminus k\subseteq\dot x_{g(\alpha)}\rbrak
  &>1-\frac{\sqrt{2\varepsilon}}{2n}\espc
  &&\text{for all $i$, for all $\alpha\in Z$,}\\
  \label{eq:59}
  \lbrak\dot y_{f_i(\beta)}\setminus k\subseteq\dot y_{h(\beta)}\rbrak
  &>1-\frac{\sqrt{2\varepsilon}}{2n}
  &&\text{for all $i$, for all $\beta\in Z$,}
\end{alignat}
where $Z$ is uncountable. Define $\dot C_i,\dot D_i\in\roone$ ($i=0,\dots,n-1$) by
$\lbrak\dot C_i(\alpha)\rbrak
=-\lbrak\dot x_{f_i(\alpha)}\setminus k\subseteq\dot x_{g(\alpha)}\rbrak$ 
and $\lbrak\dot D_i(\beta)\rbrak
=-\lbrak\dot y_{f_i(\beta)}\setminus k\subseteq\dot y_{h(\beta)}\rbrak$. 
Then by continuity and separability, 
there exists an uncountable $Z_0\subseteq  Z$ such that 
\begin{align}
  \label{eq:64}
  \lbrak n\in\dot x_{f_i(\alpha)}\rbrak&=\lbrak n\in\dot x_{f_i(\beta)}\rbrak
  &\espc\text{for all $n<k$, for all $\alpha,\beta\in Z_0$},\\
  \label{eq:74}
  \lbrak n\in\dot y_{f_i(\alpha)}\rbrak&=\lbrak n\in\dot y_{f_i(\beta)}\rbrak
  &\espc\text{for all $n<k$, for all $\alpha,\beta\in Z_0$}.
\end{align}
Let $\dot G,\dot H$ be given by $\lbrak \dot G(\alpha,\beta)\rbrak
=\lbrak\dot x_\alpha\cap\dot y_\beta\cap k\ne\emptyset\rbrak$ and 
$\lbrak\dot H(\alpha,\beta)\rbrak=\lbrak \dot x_\alpha\cap\dot y_\beta\setminus k\ne\emptyset\rbrak$. 
Thus $\dot K=\dot G\cup\dot H$, and we have 
\begin{equation}
  \label{eq:19}
  \begin{split}
  \dot K_\F&=\dot G_\F\cup\dot H_\F\\
  &\subseteq \dot G_\F\cup\bigcup_{i=0}^{n-1}\bigcup_{j=0}^{n-1}\dot K\circ (f_i,f_j)\\
  &\subseteq \dot G_\F\cup\bigcup_{i=0}^{n-1}\bigcup_{j=0}^{n-1}\dot K\circ (g,h)\cup C_i\times D_j.
  \end{split}
\end{equation}
By~\eqref{eq:64} and~\eqref{eq:74} (see~\eqref{eq:9}), 
$(Z_0,Z_0)\forces\nu(1,\dot G_\F)=0$. In conclusion,
\begin{equation}
  \label{eq:20}
  \begin{split}
  (Z_0,Z_0)\forces\nu(c,\dot K_\F)
  &\le\nu\bigl(c,\dot K\circ(g,h)\bigr)+\sum_{i=0}^{n-1}\sum_{j=0}^{n-1}\nu(c,C_i\times D_j)\\
  &<\frac\varepsilon2+n^2\cdot\frac{\varepsilon}{2n^2}\\
  &=\varepsilon.
  \end{split}
\end{equation}\\[-29pt]
\end{proof}

The following is proved in~\cite{H5}.

\begin{thm}[Hirschorn]
$\m>\aleph_1$ implies that there are no $(\oone,\oone)$ Souslin gaps in
$L^0(\nu)\div\Fin$ for a separable measure $\nu$. 
\end{thm}
 
\noindent Let us prove it here. 
Noting that $\ctbl^+\times\ctbl^+\forces b\le\lbrak m(\dot K)=1\rbrak$ (see~\eqref{eq:31}),
by Theorem~\ref{u-2},
\begin{enumeq}[(b)]
\item\label{item:4} 
$b\le\|\oone$ can be decomposed into countably many $\dot K$-homogeneous pieces$\|$.
\end{enumeq}
And by Lemma~\ref{l-16}  and Theorem~\ref{u-4},
\begin{enumeq}[(b)]
\item\label{item:5}
$-b\le\|\oone$ can be decomposed into countably 
many $\oone^2\setminus\dot K$-homogeneous pieces$\|$.
\end{enumeq}
If we make the additional assumption, as we could have, that
\begin{equation}
  \label{eq:4}
  \lbrak\dot x_\alpha\cap\dot y_\alpha\ne\emptyset\rbrak=1,
\end{equation}
then, as shown in~\cite{Hir}, by~\eqref{item:4}, $b\ne0$ implies that the 
pregap $(\dot x_\alpha,\dot y_\alpha:\alpha<\oone)$ is not Souslin; while if $b=0$
then~\eqref{item:5} implies that the pregap in $L^0(\nu)\div\Fin$ is a nongap,
completing the proof. For the reader who is at least familiar with gaps in
$\pN\div\Fin$, this simply says that $b$ forces $(\dot x_\alpha,\dot
y_\alpha:\alpha<\oone)$ is indestructible, and $-b$ forces that it is not a gap. 

\bibliographystyle{amsalpha}
\bibliography{new}

\medskip
\indent\textsc{\small Graduate School of Science and Technology,}\newline
\indent\textsc{\small Kobe University, Japan}\newline
\indent\small{\textit{E-mail address}: \texttt{James.Hirschorn@kurt.scitec.kobe-u.ac.jp}}

\end{document}